# Instability of the Environment as a Necessary Condition for Optimal Control of an Economic Object

Masaev S.N.
Siberian Federal University,
Krasnoyarsk, Russia,
Smasaev@sfu-kras.ru

*Abstract*—A review of economic approaches showed the lack of a universal method for assessing management decisions in the face of an increasing volume of analyzed data and changing parameters of the external environment. The method of integral indicators is proposed. Integral indicators are one of the modern methods for researching the behavior of an enterprise. It provides an assessment of the impact of the external environment. It shows the ability of the enterprise to adapt to new conditions. The dynamics of the correlation indicator shows the reaction of the enterprise to the impact of external factors. The purpose of the scientific work was achieved: the optimal control of the enterprise was carried out in the conditions of changing the parameters of the external environment For this, the model of the economic object and the method of its analysis are formalized. The structure of an economic object (enterprise) is given. The characteristics of the parameters of the external environment are given. The state of an economic object (enterprise) is modeled taking into account the influence of the external environment. With the help of the software package created by the author, six optimal options for control decisions have been analyzed. The state of an economic object has been modeled depending on the state of the external environment by 5,000 parameters. The research showed significant changes in the values of the correlation of the parameters of the system and the intensity of business processes when the conditions for the functioning of the system change. The optimal control of an economic object (enterprise) is selected according to the integral indicator.

## I. INTRODUCTION

Economic issues accompany the activities of any society. Without working out economic issues, it is impossible to evaluate and compare the achievements of the participants in this process: people, societies, organizations and countries. Each participant tries to optimize his activities and make it more effective by different methods [1, 2]. The activities of the participants give rise to processes: the movement of scientific and technological progress, the creation and destruction of sales markets, special economic zones [3–6], the creation and separation of holdings, the integration and conflicts of states, economic sanctions [7–10]. The result of these processes has a mutual influence on all participants in economic activity, forming random (stochastic) fluctuations and the impact of the external environment. The problem is that these fluctuations are poorly predictable or generally unpredictable. For example, the economic crises of 2008, 2014, the Covid-19 virus pandemic.

In economic theory, popular generalizations in the form of models are widespread: Walras, Sollow, circular flows, Keynesian cross, AD-AS, IS-LM, Baumol-Tobin, reproduction of Marx, Harrod-Domar, Cobb-Douglas function, multiplier-accelerator, financial accelerator, Phillips curves, Laffer curves, etc. These methods, when implemented correctly, are quite effective. Many well-known specialists have dealt with the issues of modeling and control economic systems, starting with the Nobel laureates V.V. Leontyev and L.V. Kantorovich, and ending with modern economists such as A.G. Granberg, A.G. Aganbegyan, V.F. Krotov and others [11, 12]. These scientific works are classical and define the object of study. Formal statement of the problem of observation and control.

The above described tasks do not cause difficulties for solving. However, there are economic entities that are characterized by a high dimension of big data. On them, the simplest economic task becomes practically impossible. For example, business plans from five hundred to several thousand pages.

These difficulties complicate the work of specialists making decisions about the assessment of an investment project or managers making decisions.

From the above, it can be argued that there is a complication. As a result, world leaders companies have more processes in progress than some of the leading countries or countries in the early 20th century. The creation of a general indicator that takes into account all indicators of an economic object and external changes is relevant.

The growing dimension of economic objects forms an urgent research purpose: to perform optimal control of an enterprise in conditions of changing parameters of the external environment.

The concept of "optimal" is understood as a control decision that minimizes/maximizes the final indicator on a known and adjustable set of parameters.

The achievement of the research set in the article depends on the tasks:

- Formalization of the economic object model and the method of its analysis;

- Description of the structure of the economic object;

- Description of the parameters of the external environment;

- Modeling and analysis of the state of an economic object.

In 2008, the method of integral indicators was created. The method of integral indicators was tested in the crisis of 2008. This method was used to predict the world financial crisis. The integral method captures positive or negative phenomena in the observed object. The method was widely used: for the analysis of actual data on the activities of banks and stock exchanges [13], insurance markets [14], investment planning [15], the creation of special economic zones [16], training of employees of the penal system [17], performance evaluation business processes [18], search for information anomalies [19], training of employees in departmental universities [20], etc. [21, 22].

In 1987 A.N. Gorban and his co-workers proposed the method of correlation adaptometry. The method made it possible to assess the adaptation of living organisms to stress. Stress means both good and bad situations. The method of correlation adaptometry is the basis for the method of integral indicators.

We will use the method of integral indicators to analyze an economic object, so it makes sense to consider it briefly when formalizing an economic enterprise (enterprise).

## II. RESEARCH METHODS AND MATERIALS

For the research, the activity of the enterprise is considered as a system of business processes $S$

$$S=\{T,X\},$$

where $T=\{t: t=1,\ldots, T_{max}\}$ – many points in time (months);

$X$ – system parameter space;

$x(t)=[x^1(t), x^2(t),\ldots, x^n(t)]^T \in X$- $n$ – system state vector;

$x_i(t)$ – enterprise income/expenses.

The dimension $n=5,000$ parameters.

To analyze the system at time $t$, we consider the parameters $x(t)$ for the previous clock cycles $k$.

Then

$$X_k(t) = \begin{bmatrix} x_1^1(t-1) & x_2^1(t-1) & \ldots & x_n^1(t-1) \\ x_1^2(t-2) & x_2^2(t-2) & \ldots & x_n^2(t-2) \\ \ldots & \ldots & \ldots & \ldots \\ x_1^i(t-k) & x_2^i(t-k) & \ldots & x_n^i(t-k) \end{bmatrix} \quad (1)$$

The correlation matrix is calculated for the maximum analysis horizon

$$R_k(t) = \frac{1}{k-1} \overset{o}{X_k^T}(t) \overset{o}{X_k}(t) = \|r_{ij}(t)\| \quad i,j=1,\ldots,n, \ (2)$$

where

$$\overset{o}{r_{ij}}(t) = \frac{1}{k-1} \sum_{l=1}^{k} \overset{o}{x^i}(t-l) \overset{o}{x^j}(t-l) \quad (3)$$

$t$ – moments in time; $r_{ij}(t)$ – correlation coefficients of variables $x^i(t)$ and $x^j(t)$ at the moment $t$.

Following the method of integral indicators, we will form the indicator of the correlation coefficients into an indicator of express assessment of the correlation of the parameters of the economic system $G_i(t)$:

$$R_i(t) = G_i(t) = \sum_{j=1}^{n} |r_{ij}(t)| : (|r_{ij}(t)| \ge r_{sign}), \quad (4)$$

where $r_{sing}$ the critical value of the correlation coefficient at a given depth of analysis $k$.

An integral indicator of the entire system of the $i$-th enterprise.

$$G = \sum_{t=1}^{T=max} \sum_{i=1}^{n} G_i(t) \quad (5)$$

## III. DESCRIPTION OF THE RESEARCH OBJECT (STRUCTURE OF THE ECONOMIC OBJECT)

The modeling was carried out at the enterprise for deep processing of wood. The company operates in the Severo-Yenisei region. The main business processes of the enterprise: harvesting 800 thousand cubic meters of round timber, delivery of round timber along the Yenisei River, production of deep processing products from round timber. All business processes are described according to the structure of the business plan [23].

The business plan has 1060 pages. It is a priority project of the Russian Federation (order on the second page of the presentation [23]).

The examination of the business plan was signed by the Minister of Economy of the Russian Federation. Published by the Trade Mission of the Russian Federation in China [23].

Due to the large amount of data, the descriptive part of the business plan, we will only give the structure of the formation of the relationship between the parameters of the enterprise. The zero block is the parameters of the external environment: exchange rates, taxes of the Russian Federation, the cost of fuel and electricity tariff, the size of inflation for resources and products, prices for raw materials and products, prices for equipment, the level of equipment productivity.

The first block is an investment plan and a listing of all ongoing activities and business processes at the enterprise.

The second block is equipment, machines and mechanisms used in production (harvesters, forwarders, tractors, excavators, graders, timber trucks "chance", diesel power generators, rotational cars, buildings and structures, equipment for deep processing of USNR wood), machinery and org. technique used for administrative staff. All types of equipment repair have been modeled: major, current and operational inspection, depending on the mileage of each vehicle or its use by administrative personnel. Equipment maintenance is based on production over five years.

The third block – depreciation is calculated for all fixed assets and equipment.

The fourth block is a description of the characteristics and volumes of manufactured products: round timber, calibrated timber, glued board, euro lining, floor boards, veneer, furniture and pellets.

The fifth block is the management of balances in the warehouse and the logistics of moving material assets, taking into account the main technological stages: a) procurement, b) storage at the upper section, c) delivery by barges in May, June, July, August, September along the Yenisei River to the lower warehouse , d) storage at the lower platform, e) deep processing of round timber and production of: calibrated timber, glued board, euro lining, floor boards, veneer, furniture) pellet production - waste processing.

The sixth block is the formation of the regular number of workers, engineering and technical personnel, administrative personnel and other employees.

The seventh block is the calculation of the budget and the economy, loans, loans, the owner's funds at the enterprise, the formation of Form 2 (Profit and Loss Statement), cash flow and Form 1 (Balance Sheet). The model calculates all existing economic indicators to evaluate the enterprise. The eighth block is the calculation of environmentally friendly production parameters. The ninth block is the engineering of the entire production by USNR (an American company, the world leader in the production of machinery and equipment for deep processing of round timber).

In the research, the enterprise model changes its state at each time period depending on environmental conditions.

## IV. ENVIRONMENTAL PARAMETERS

Parameters (factors) of the external environment: sanctions on high-tech equipment, a ban on the sale of round timber for export, government subsidies, the size of a loan from a bank, a bank rate on a loan, Sanctions on the supply of high-tech equipment shift the implementation of the project by 5 months, as time is spent searching similar equipment. The introduced ban on the export of round timber in Russia leads to an increase in the project. Banks are forced to raise interest rates for loans due to new risks. The enterprise requires more credit resources, more funds from the owner, assistance from the state is needed in the form of subsidies, the degree of processing of round timber for export sales.

The interaction of these factors is combined into six options for managing enterprise development. Control options are presented below.

First option. Attracted banker confirmed credit of the project cost 100% at basic lending rate 10% per annum. There is no owner investment. Subsidies from the state are 33% of the amount of subsidies required by law. The sale of sawlogs starts from the 5th month from the start of the project. Sale of deep processing products in 21.

Second option. Attracted banker confirmed credit of the project cost 63% at basic lending rate 10% per annum. Investments of the project owner 37% of the project cost. Subsidies from the state are 100% of the amount of subsidies required by law. The sale of sawlogs starts from the 10th month from the start of the project. Sale of deep processing products from 21 months from the start of the project.

Third option. Attracted banker confirmed credit of the project cost 63% at basic lending rate 10% per annum. Investments of the project owner 37% of the project cost. Subsidies from the state are 100% of the amount of subsidies required by law. The sale of sawlogs starts from the 10th month from the start of the project. Sale of deep processing products from 21 months from the start of the project. The offset of the purchase of fixed assets from the original date is 5 months.

The fourth option. Attracted banker confirmed credit of the project cost 126% at basic lending rate 13% per annum. There is no owner investment. Subsidies from the state are 100% of the amount of subsidies required by law. The sale of sawlogs starts from the 10th month from the start of the project. Sale of deep processing products from the 27th month from the start of the project. The offset of the purchase of fixed assets from the original date is 5 months.

Fifth option. Attracted banker confirmed credit of the project cost 112% at basic lending rate 13% per annum. There is no owner investment. Subsidies from the state are 100% of the amount of subsidies required by law. The logs are shipped for sale from the tenth month from the start of the project. Products have been sold since the twenty-seventh month. Fixed assets are purchased from the fifth month of the project.

Sixth option. Attracted banker confirmed credit of the project cost 112% at basic lending rate 13% per annum. There is no owner investment. Subsidies from the state are 100% of the amount of subsidies required by law. The logs are shipped for sale from the tenth month from the start of the project. Products have been sold since the twenty-seventh month. Fixed assets are purchased from the fifth month of the project. Increase in logging from 800 to 1,000 thousand cubic meters

After analyzing the parameters of production options, you can determine the impact of the external environment on each parameter.

The amount of credit for the implementation of production is one of the first parameters that respond to the deterioration of environmental conditions and changes from the second to the sixth periods. The annual interest rate on the loan changes only from the fourth period to the sixth period.

The ban on the sale of roundwood affects the second to sixth options.

The imposed sanctions on high-tech equipment in the United States affect the third through sixth options.

Obviously, in the first variant, environmental conditions do not cause significant changes in the parameters of production. In subsequent implementations, environmental conditions increasingly affect production.

The issue of sanctions is very complex. It is considered in a separate work [24].

## V. FINDINGS

Modeling and calculations were carried out for 5,000 parameters using the software package developed by the author.

The Figure 1 shows dynamics integral indicator

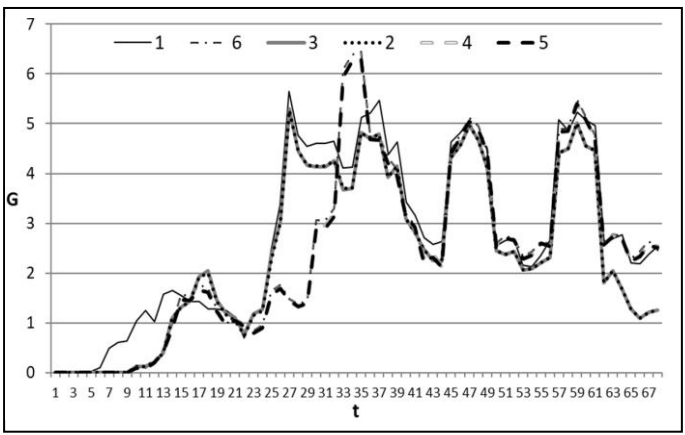

Fig. 1. Dynamics $G$

When calculating a large number of relationships, enterprise analysis is difficult. It is expedient to display a large number of correlations by the graph $G_i(t)$. This is done in the program.

From 5 to 21 periods this is the preparation of production for the release of products. From 25 to 37 periods directly the start of production. Changes in the dynamics of the integral indicator from 45 to 50, from 57 to 62 periods are the standard operating mode of the enterprise with the influence of the factors of the external environment.

Each point of the graph in Figure 1 can be displayed as a correlation graph. Figure 2 shows significant mutual correlations of elements $x^i(t)$ and $x^j(t)$ in the 18th period for the 2nd control option. Figure 3 shows the correlation graph in the 35th period for the same control option, respectively. Figures 2 and 3 show the normal operation of the enterprise and during the period of adaptation to external environmental parameters.

It is characterizing the number of relationships in the correlation graph.

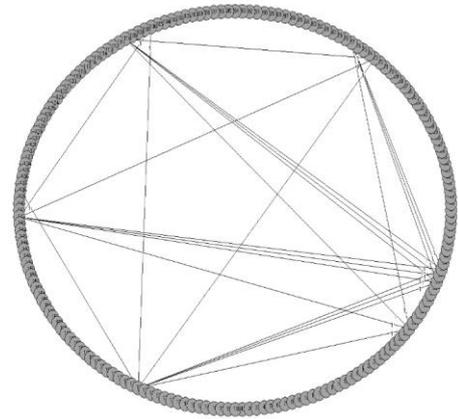

Fig. 2. Correlation graphs for $t$=18

A Figure 3 shows the correlation graph

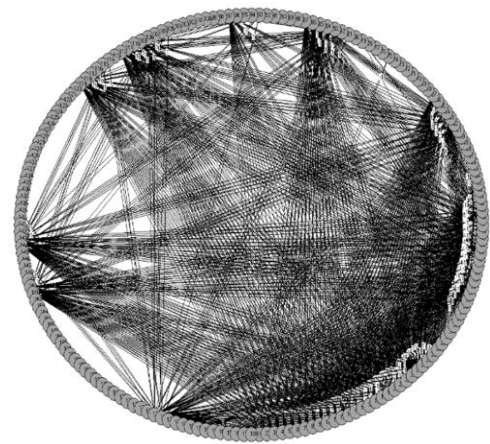

Fig. 3. Correlation graphs for $t$=35

There are three ways from the polyfactorial condition, assessed by A. Gorban and E. Smirnova: 1 – the correlation increases with an increase in the influence of the external environment, 2 – adjustment to new factors and further work, 3 – non-adjustment to external factors and death. The growth of $G$ characterizes the significant influence of environmental factors. Table 1 shows a decrease in $G$ characterizes the adaptation processes to new environmental conditions.

TABLE I. VALUES OF THE INTEGRAL INDICATOR

| Parameters | Option | $G$ |
|---|---|---|
| Integral indicator | First | 186.6 |
| | Second | 161.7 |
| | Third | 162.0 |
| | Fourth | 162.8 |
| | Fifth | 162.5 |
| | Sixth | 166.5 |

With our company as well. With an increase in the influence of external factors on its business processes, the value of the integral indicator grows. Figure 2and 3 show it. Each point characterizes one business process of the enterprise. Figure 2 shows the state of the business processes of the enterprise under the weak influence of the external environment. Figure 3 shows the state of business processes under the strong influence of the external environment.

Figure 1 shows that with the implementation of all production options of the enterprise, favorable and unfavorable situations arise. According to the integral indicator, the first variant of production implementation is preferable. From the 32nd period to the 38th period in the fourth, fifth and sixth variants, crisis phenomena were recorded as an integral indicator. The integral indicator records the imposition of sanctions on equipment purchased from the United States.

The index reflects the influence of seasonal factors from the 43rd to the 50th and from the 56th to the 62nd periods: inflation, delivery of raw materials along the Yenisei River.

## VI. DISCUSSION

The dynamics of the integral indicator characterizes the activity of economic objects. Processes can be evaluated:

- Bankruptcy of economic objects;
- Establishment of enterprises;
- Integrations, acquisitions, etc.

It is shown that the size and type of activity of the analyzed object does not matter for the integral indicator. Therefore, it can be applied to analyze various economic objects.

The dynamics of the integral indicator reflects well the influence of the external environment on the activity of an economic object. Therefore, it is possible to analyze plans for the development of an economic object in different situations.

The external environment is changing so quickly that recurring factors began to appear that change the activity of an economic object. For example, the situation with the COVID-19 pandemic. Countries and large companies were not ready for it. However, the second wave of the pandemic does not cause such horror.

It is also worth noting that in addition to changes in the external environment, it is worth considering changes in internal factors. This change can also lead to significant deviations in the activities of an economic entity.

An integral indicator can also be used to track internal changes in an economic object. For example, the choice of a management method and its assessment:

- The Interstate Standardization System;
- Soviet quality management system;
- US National Aeronautics and Space Administration Regulations;
- The Yitzhak Adizes method;
- The acquired personnel of the RF national qualifications framework;
- Dublin Descriptors;
- Cognitive abilities of staff (B. Bloom's taxonomy);
- Acquired universal competencies of personnel;

- Competencies of employed graduates;
- Strategic planning;
- Job descriptions;
- BSi quality method personnel competence (UK);
- Goal classifier;
- European Qualifications Framework;
- British Institute for Quality Method;
- Competence of the personnel of the project method Hermes (Switzerland);
- Project method PROMAT (South Korea);
- P2M project method (Japan);
- Design method ANCSPM (Australia);
- Models of information systems development V-Modell, Germany);
- Design method DIN 69901 (Germany);
- Project method PRINCE2 (UK);
- Personnel competencies of the P2M project method (Japan);
- VZPM project method (Switzerland);
- R. Solow models;
- Competencies of the personnel of the PMBoK design method (USA);
- Competence of the personnel of the IPMA design method (Switzerland);
- Design method CAN / CSA-ISO 10006-98 (Canada);
- Scripting method Hermes (Switzerland).

## VII. CONCLUSION

The research provides a general assessment of all business processes of an enterprise with an integral indicator. A model of a wood processing enterprise is used as an example. There are six options for the development of the enterprise. Each variant is evaluated by an integral indicator.

It was found that the integral indicator assesses the impact of environmental parameters on each option for the implementation of production in different ways. The first variant of production implementation is the priority for application.

The integral indicator characterizes all indicators of the enterprise. The model of the enterprise has passed expertise in the federal ministries. Then we can assume that the proposed integral indicator can be applied to other enterprises and various areas of analysis.

So, the tasks of the article: formalization of the economic object model and the method of its analysis, description of the structure of the economic object, description of the parameters of the external environment, modeling of the economic object, have been achieved. Consequently, the purpose of the research: to perform optimal enterprise management in conditions of changes in the parameters of the external environment, has been achieved.